\documentclass[12pt,a4paper]{article}
\usepackage[utf8]{inputenc}
\usepackage{amsmath}
\usepackage{bm}
\usepackage{amsfonts}
\usepackage{amssymb}
\usepackage{graphicx}
\usepackage[colorlinks,citecolor = brown,linkcolor=violet,anchorcolor = blue]{hyperref}
\usepackage{tikz}
\usepackage{parskip}
\usepackage{amsthm}
\usepackage[toc,page,titletoc]{appendix}

\usepackage{wrapfig}
\usepackage{mathtools}
\usepackage{graphicx}
\graphicspath{{figures/}} 
\usepackage{subcaption}
\usepackage{array}
\usepackage{dirtytalk}
\usepackage{arydshln} 
\numberwithin{equation}{section}
\usepackage{cleveref}
\usepackage[sort&compress,numbers]{natbib}
\usepackage{lipsum}
\usepackage{placeins}

\newcommand{\Order}[1]{\mathcal{O}(#1)}
\usepackage{xcolor}

\usepackage{xparse} 
\ExplSyntaxOn
\NewDocumentCommand{\longdash}{ O{2} }
 {
  --\prg_replicate:nn { #1 - 1 } { \negthinspace -- }
 }
\ExplSyntaxOff

\crefname{supp}{Supplement}{Supplements}

\crefformat{section}{\S#2#1#3} 
\crefformat{subsection}{\S#2#1#3}
\crefformat{subsubsection}{\S#2#1#3}

\usepackage{geometry}
\author{Sarafa A. Iyaniwura\textsuperscript{1} and Zhiwei Peng\textsuperscript{2}}

\title{Asymptotic analysis and simulation of mean first passage time for active Brownian particles in 1-D}

\date{}

\begin{document}

\maketitle
\small{
\textbf{1} Department of Mathematics and Institute of Applied Mathematics, University of British Columbia, Vancouver, BC,  V6T 1Z2, Canada.\\
\textbf{2} Department of Chemistry, University of Toronto, Toronto, ON, M5S 3H6,
Canada.
}

*\textbf{Correspondence:} iyaniwura@aims.ac.za \& zpeng@alumni.caltech.edu

\begin{abstract}
Active Brownian particles (ABPs) are a model for nonequilibrium systems in which the constituent particles are self-propelled in addition to their Brownian motion. Compared to the well-studied mean first passage time (MFPT) of passive Brownian particles, the MFPT of ABPs is much less developed.  In this paper, we study the MFPT for ABPs in a 1-D domain with absorbing boundary conditions at both ends of the domain. To reveal the effect of swimming on the MFPT, we consider an asymptotic analysis in the weak-swimming or small P\'eclet ($Pe$) number limit. In particular, analytical expressions for the survival probability and the MFPT are developed up to $\mathcal{O}(Pe^2)$. We explore the effects of the starting positions and starting orientations on the MFPT. Our analysis shows that if the starting orientations are biased towards one side of the domain, the MFPT as a function of the starting position becomes asymmetric about the center of the domain. The analytical results were confirmed by the numerical solutions of the full PDE model.
\end{abstract}

\textbf{Keywords}: Mean first passage time, asymptotic analysis, active Brownian particles, survival probability.

\section{Introduction}
Mean first passage time (MFPT) is the average time-scale for a stochastic event to first occur \cite{van1992stochastic, redner2001guide}.
This phenomenon has been applied to many physical and biological problems \cite{fauchald2003using, mckenzie2009first, godec2016first, metzler2014first, polizzi2016mean}, such as calculating the time it takes for a protein to find a binding site on DNA \cite{mirny2009protein}, the time it takes for a predator to find its prey \cite{kurella2015asymptotic}, and computing the time it takes for a diffusing molecule to reach a localized  signaling region on a cell membrane \cite{coombs2009diffusion}, among others. 
The MFPT for Brownian particles has been studied extensively (see \cite{pillay2010asymptotic, condamin2005first, lindsay2017optimization, cheviakov2010asymptotic, iyaniwura2021optimization, mattos2012first,
lindsay2017first, iyaniwura2021simulation},  and  the  references  therein). However, not much has been done on the MFPT for active Brownian particles (ABPs). An active Brownian particle is a model of the self-propelled motion of active matter systems, such as motile bacteria and synthetic active particles \cite{romanczuk2012active,Bechinger2016}. In addition to the translational Brownian motion of passive particles, an ABP exhibits rotational Brownian motion and self-propulsive (or swimming) motion. At long times, the swimming direction (or orientation) of ABPs is randomized due to rotational Brownian motion. A second model for bacteria locomotion that is closely related to the ABP model is the so-called run-and-tumble particle (RTP) model. Instead of continuous rotational diffusion, RTPs undergo discrete tumbling events  that randomize their swimming direction \cite{solon2015active}. Some articles that looked at the first passage time for ABPs include \cite{moen2022trapping, scacchi2018mean, angelani2015run, bressloff2023trapping}. In \cite{scacchi2018mean}, the MFPT for ABP  was studied numerically in a 1-D domain. If the particle's orientation is restricted to 1-D, the orientation space becomes binary, i.e., the particle can either take the right orientation or the left orientation. Due to the discrete orientation space in 1-D, no distinction will be made between the ABP and RTP models. For both models, the random orientational dynamics reduces to a random walk between positive and negative orientations.

The first passage time distribution and its mean have been calculated for passive Brownian particles in different domains \cite{lindsay2017first, grebenkov2016universal, cheviakov2010asymptotic, iyaniwura2021optimization, condamin2007first, pillay2010asymptotic, tzou2015mean, grebenkov2020mean, mattos2012first, iyaniwura2021simulation}. Many of these studies used asymptotic analysis to calculate the MFPT for a Brownian particle to escape a domain, either through a small absorbing trap localized in the domain or through a small opening on the boundary of the domain \cite{cheviakov2010asymptotic, lindsay2017first, pillay2010asymptotic}. Other studies have looked at optimizing trap configurations that minimize MFPT in different geometries for  stationary  and moving traps \cite{iyaniwura2021optimization, iyaniwura2021simulation, tzou2015mean, kolokolnikov2005optimizing}. In contrast to the extensive research on the FPT/MFPT for Brownian particles, the FPT for active Brownian particles \cite{moen2022trapping, scacchi2018mean, angelani2015run} is much less understood. Notably, analytical progress is lacking in part due to the added complexity of rotational Brownian motion and swimming. For passive Brownian motion, one can often start the analysis from the partial differential equation governing the FPT distribution \cite{redner2001guide}. For active Brownian motion, such an equation is not yet available, and one needs to consider the full probability density function of finding the ABP at a given position, orientation, and time.  With a solution to the probability density function in space and time, one can then obtain the survival probability and the resulting MFPT.  

In this study, we consider the MFPT of ABPs in a 1-D domain with absorbing boundary conditions at both ends of the domain. The main aim of this study is to understand the effect of the swimming speed of ABPs on the MFPT. If the swim speed is zero, one recovers the MFPT of passive Brownian particles. To understand how the swimming motion affects the MFPT, we first consider an  asymptotic analysis in the weak-swimming limit. This analysis allows us to reveal the first effects of swimming on the MFPT distribution in contrast to the results for passive Brownian particles.

Consider a 1-D domain $\Omega \equiv [-R, R]$, where $R \in \mathcal{Z}^+$, and let $P_+(x,t)$ and $P_{-}(x,t)$ be the  density of positive-oriented (pointing to the right direction) and negative-oriented (pointing to the left direction) active Brownian particles (ABPs), respectively, at position $x$ at time $t$. 
Due to the run-and-tumble nature of ABPs, positive-oriented particles can change orientation to the left. Likewise, the negative-oriented particles can change their orientation to the right. Based on this, the dynamics of  $P_+$ and $P_{-}$ satisfy the following partial differential equations (PDEs)
\begin{subequations}\label{Eq:Dim_Model_P}
\begin{align}
\frac{\partial P_+}{ \partial t}\,& = -\frac{\partial}{ \partial x} \Big{(} v_s\,P_+  - D_T\, \frac{\partial P_+}{\partial x}\Big{)} - \frac{1}{\tau} P_+ + \frac{1}{\tau} P_{-},  \quad x \in \Omega, \quad t > 0;\label{Eq:1D_ABP_1}\\[2ex]
 \frac{\partial P_{-}}{ \partial t}\,& = -\frac{\partial}{ \partial x} \Big{(} -v_s\,P_{-}  - D_T\, \frac{\partial P_{-}}{\partial x}\Big{)} + \frac{1}{\tau} P_{+} - \frac{1}{\tau} P_{-},\label{Eq:1D_ABP_2}
\end{align}
where $v_s$ and $D_T$ are the dimensional uniform swimming speed and diffusivity of the ABPs, respectively, and 
$1/\tau$ is the tumbling rate of the particles, that is, the rate at which the particles switch from positive orientation to negative orientation and vice versa. This implies that $\tau$ is the time it takes for a particle to change its orientation. In \eqref{Eq:1D_ABP_1} and \eqref{Eq:1D_ABP_2}, we have assumed that the particles change their orientations from right to left and vice versa at the same rate.
We specify Dirichlet boundary conditions, $P_{\pm}(\pm R, t) = 0$, at the two ends of the domain. These boundary conditions impose that an ABP `vanishes' or is absorbed when it hits the boundaries. In addition, we impose the following initial conditions
\begin{equation}
\label{Eq:initial-condition-P}
    P_{+}(x, 0) = \eta \,\delta(x-x_0) \quad \text{and}  \quad P_{-}(x, 0) = (1 - \eta)\, \delta(x-x_0),
\end{equation}
\end{subequations}
where $\delta(x-x_0)$ is the Dirac delta function localized at $x_0$, and $\eta \in [0,1]$ is the fraction of particles located at $x_0$ at $t=0$ with positive-orientation. 
This model follows a similar framework to the model of \cite{angelani2015run, angelani2014first}. However, their model only considers run-and-tumble particles without diffusion ($D_T=0$). We define the total density of particles at position $x$ at time $t$
as $n(x,t) = P_+(x,t) + P_{-}(x,t)$ and $f(x,t) = P_+(x,t) - P_{-}(x,t)$.
Upon adding \eqref{Eq:1D_ABP_1} and \eqref{Eq:1D_ABP_2}, and subtracting \eqref{Eq:1D_ABP_2} from \eqref{Eq:1D_ABP_1} in two different operations, we obtained a coupled PDE system for $n(x,t)$ and $f(x,t)$, given by 
\begin{subequations}\label{Eq:Dim_Model_n}
\begin{align}
    \frac{\partial n}{ \partial t}\,&= -\frac{\partial}{ \partial x} \Big{(} v_s\,f  - D_T\, \frac{\partial n}{\partial x}\Big{)}, \quad x \in \Omega, \quad t > 0 \, ; \label{Eq:N_eqn}\\[2ex]
        \frac{\partial f}{ \partial t}\,&= -\frac{\partial}{ \partial x} \Big{(} v_s\,n  - D_T\, \frac{\partial f}{\partial x}\Big{)} - \frac{2}{\tau}f. \label{Eq:f_eqn}    
\end{align}
As shown in \eqref{Eq:N_eqn}, active Brownian particles are transported by their swimming motion in addition to Brownian diffusion. One can recover the diffusion equation for the total density of particles $n(x,t)$ for passive Brownian particles by setting $v_s =0$ in \eqref{Eq:N_eqn}. 
The boundary and initial conditions for the coupled PDE system \eqref{Eq:N_eqn} and \eqref{Eq:f_eqn} are given by
\begin{align}
   n(\pm R, t) = 0 & \quad \text{and} \quad f(\pm R, t) = 0; \label{Eq:Dim_Model_BC} \\
 n(x, 0) = \delta(x-x_0) \quad & \text{and} \quad f(x, 0) = (2\eta-1) \, \delta(x-x_0). \label{Eq:Dim_Model_IC}
\end{align}
\end{subequations}
Our goal is to use \eqref{Eq:Dim_Model_n} to study the MFPT for active Brownian particles in $\Omega$, in the weak-swimming limit.
We define the survival probability, $S(t; x_0)$ as the probability that an ABP that starts at position $x_0 \in \Omega$ at time $t=0$ is still in the domain $\Omega$ at time $t > 0$. It is given by
\begin{equation}\label{Eq:Surv_Prob}
    S(t; x_0) = \int_{\Omega} n(x,t) \;\text{d}x.
\end{equation}
where $n(x,t)$ is the total density of particles at position $x$ at time $t$.
In terms of the survival probability $S(t; x_0)$, we define the first passage time (FPT) distribution $F(t; x_0)$ for an ABP, starting  at position $x_0 \in \Omega$ at time $t=0$, to escape the domain $\Omega$ (through its boundaries) as 
\begin{align}\label{Eq:fpt-distribution}
    F(t; x_0) = - \frac{\partial S(t; x_0)}{\partial t}.
\end{align}
Similarly, we define the MFPT for an ABP as the average time scale for the particle to escape the domain, starting from a position within the domain. For a particle starting at position $x_0 \in \Omega$ at time $t=0$, the MFPT, $\mathcal{\mu}(x_0)$, is computed from the first passage time distribution \eqref{Eq:fpt-distribution} as follows \cite{redner2001guide} 
\begin{align}\label{Eq:MFPT_formula}
  \mathcal{\mu}(x_0) =  \int_0^\infty \tau F(\tau; x_0) \, \text{d} \tau.  
\end{align}
Using \eqref{Eq:Surv_Prob} and \eqref{Eq:fpt-distribution}, we derive an expression for the MFPT, $\mu(x_0)$ in terms of the survival probability, $S(t; x_0)$, given by
\begin{align}\label{Eq:MFPT_formula_Surv}
  \mu(x_0) = \int_0^\infty S(t; x_0) \; \text{d} t,
\end{align}
provided that $S(t)$ decays to zero faster than $1/t$ as $t \to \infty$. This formula will be used to compute the MFPT in this study.

\section{Weak-swimming asymptotic analysis}\label{sec:1D_WeakSwimAsymp}

We analyze the coupled PDE system \eqref{Eq:Dim_Model_n} in the weak-swimming regime. In this regime, we assume that the swimming speed of the active Brownian particles is small and use asymptotic analysis to derive an approximate solution to the coupled system. We non-dimensionalize the PDEs in \eqref{Eq:Dim_Model_n} by scaling $x$ with the length-scale of the domain $R$ and $t$ with the diffusive time-scale $R^2/D_T$ to obtain
\begin{subequations}\label{Eq:Dimless_Model_n}
\begin{align}
    \frac{\partial n}{ \partial t}\,&= -\frac{\partial}{ \partial x} \Big{(} Pe\,f  -  \frac{\partial n}{\partial x}\Big{)}, \quad x \in \Omega_{u} \equiv [-1,1], \quad t > 0\,; \label{Eq:Dimless_N_eqn}\\[1ex]
     \frac{\partial f}{ \partial t}\,&= -\frac{\partial}{ \partial x} \Big{(} Pe\,n  - \frac{\partial f}{\partial x}\Big{)} - 2\beta f \,; \label{Eq:Dimless_f_eqn}\\[1ex]
     n(\pm &1, t)= 0 \quad \text{and} \quad f(\pm 1, t) = 0 \,; \label{Eq:Dimless_BC}\\[1ex]
     n(x,0) & = \delta(x-x_0) \quad  \text{and} \quad f(x, 0) = (2\eta-1) \, \delta(x-x_0), \label{Eq:Dimless_IC}
\end{align}
\end{subequations}
where $Pe = (v_s R)/D_T$ is the swimming P\'eclet number and $\beta= R^2 /(\tau D_T)$. For fixed values of $R$ and $D_T$, $Pe$ is directly proportional to the dimensional swimming speed, and $\beta$ is proportional to the tumbling rate of the particles. We  analyze the dimensionless coupled PDE system \eqref{Eq:Dimless_Model_n} in the weak-swimming limit, i.e., for $Pe \ll 1$, using asymptotic analysis.

We begin our analysis by expanding $n(x,t)$ and $f(x,t)$ in terms of $Pe \ll 1$ as follows
\begin{subequations}\label{Eq:Dimless_Expand}
\begin{align}
    n(x,t)\,&= n_0(x,t) + Pe\, n_1(x,t) + Pe^2 \, n_2(x,t) +  \cdots \label{Eq:Dimless_N_expand}\\[1ex]
    f(x,t)\,&= f_0(x,t) + Pe\, f_1(x,t) + Pe^2 \, f_2(x,t) + \cdots \label{Eq:Dimless_f_expand}
\end{align}
\end{subequations}
Similarly, we expand the survival probability $S(t)$ and MFPT $\mu$ as 
\begin{subequations}\label{Eq:Surv_FPT_Expand}
\begin{align}
     S(t) &= S_0(t) + Pe\, S_1(t) + Pe^2\, S_2(t) + \cdots  \label{Eq:SurvivalExpand}\\[1ex]
    \mu &= \mu_0 +Pe \,\mu_1 + Pe^2 \,\mu_2 + \cdots \label{Eq:MFPT_expand}
\end{align}
\end{subequations}
Note that the survival probability ($S(t)$) and the MFPT ($\mu$) are dimensionless since $n$ and $f$
 are dimensionless. In particular, $\mu$ is obtained by scaling the dimensional MFPT by the diffusive time scale $R^2/D_T$. 
Upon substituting \eqref{Eq:Dimless_Expand} into \eqref{Eq:Dimless_Model_n} and collecting terms in powers of $Pe$, we obtain the leading-order problem, given by
\begin{subequations}\label{Eq:WkSwim_LD_Prob}
\begin{align}
    \frac{\partial n_0}{ \partial t}\,&=  \frac{\partial^2 n_0}{\partial x^2},  \quad x \in \Omega_{u}, \quad t > 0\,;\label{Eq:WeakSwim_leadOrderN}\\[1ex]
        \frac{\partial f_0}{ \partial t}\,&=  \frac{\partial^2 f_0}{\partial x^2} - 2\beta\,f_0 \,;\label{Eq:WeakSwim_leadOrderF}\\[1ex]
           n_0(\pm 1, t)&= 0  \quad \text{and} \quad f_0(\pm 1, t) = 0 \,; \label{Eq:WkSwim_LD_BC}\\[1ex]
n_0(x, 0) =&  \delta(x-x_0) \quad  \text{and} \quad f_0(x, 0) = (2\eta-1) \, \delta(x-x_0)\label{Eq:WkSwim_LD_IC}, 
\end{align}
\end{subequations}
We notice from \eqref{Eq:WkSwim_LD_Prob} that the leading-order PDEs are decoupled, and the density of particles $n_0(x,t)$ satisfies the diffusion equation. 
At $\mathcal{O}(Pe)$, we have
\begin{subequations}\label{Eq:WkSwim_PesOrder_Prob}
\begin{align}
    \frac{\partial n_1}{ \partial t}\,&=  \frac{\partial^2 n_1}{\partial x^2} - \frac{\partial f_0}{\partial x},  \quad x \in \Omega_{u} \equiv [-1,1], \quad t > 0\,;\label{Eq:WeakSwim_Pes_Order_n} \\[1ex]
    \frac{\partial f_1}{ \partial t}\,&= - \frac{\partial n_0}{\partial x} +   \frac{\partial^2 f_1}{\partial x^2} - 2\beta\,f_1\,; \label{Eq:WeakSwim_Pes_Order_f}\\[1ex]
     n_1(\pm 1,&\, t)= 0  \quad \text{and} \quad f_1(\pm 1, t) = 0 \,;\label{Eq:WkSwim_Pes_Order_BC}\\[1ex]
n_1(x,&\, 0)  =  0 \quad  \text{and} \quad f_1(x, 0) = 0\label{Eq:WkSwim_Pes_Order_IC}.
\end{align}
\end{subequations}
The $\mathcal{O}(Pe^2)$ problem is given by 
\begin{subequations}\label{Eq:WkSwim_Pes2Order_Prob}
\begin{align}
    \frac{\partial n_2}{ \partial t}\,&=  \frac{\partial^2 n_2}{\partial x^2} - \frac{\partial f_1}{\partial x},  \quad x \in \Omega_{u} \equiv [-1,1], \quad t > 0 \,; \label{Eq:WeakSwim_Pes2_Order_n} \\[1ex]
    \frac{\partial f_2}{ \partial t}\,&= - \frac{\partial n_1}{\partial x} +   \frac{\partial^2 f_2}{\partial x^2} - 2\beta\,f_2\,; \label{Eq:WeakSwim_Pes2_Order_f}\\[1ex]
     n_2(\pm 1,&\, t)= 0  \quad \text{and} \quad f_2(\pm 1, t) = 0 \,;\label{Eq:WkSwim_Pes2_Order_BC}\\[1ex]
n_2(x,&\, 0)  =  0 \quad  \text{and} \quad f_2(x, 0) = 0. \label{Eq:WkSwim_Pes2_Order_IC}
\end{align}
\end{subequations}
We shall solve each of the problems \eqref{Eq:WkSwim_LD_Prob}, \eqref{Eq:WkSwim_PesOrder_Prob}, and \eqref{Eq:WkSwim_Pes2Order_Prob}, and use their solutions to construct a three-term asymptotic expansion for the MFPT.

We begin with the leading-order problem \eqref{Eq:WkSwim_LD_Prob}.
Using separation of variables, we obtain that $n_0(x,t)$ satisfies  
\begin{equation}\label{Eq:WkSwim_n0_Sol}
     n_0(x,t) = \sum_{n=0}^{\infty} \Big{[}  \cos \left( \lambda_{1,n} \, x_0 \right) \cos \left( \lambda_{1,n} \, x \right)e^{-\lambda_{1,n}^2 \, t } + \sin(\lambda_{2,n} x_0) \sin \left( \lambda_{2,n} x \right) e^{-\lambda_{2,n}^2 \, t } \Big{]}.
\end{equation}
where $\lambda_{1,n} = (2n + 1)\pi/2$ and $\lambda_{2,n} = n\pi$ for $n \in \mathbb{Z}$.
We use the survival probability formula in \eqref{Eq:Surv_Prob} together with the asymptotic expansion in \eqref{Eq:SurvivalExpand}  to construct the leading-order survival probability, $S_0(t; x_0)$ as an integral of $n_0(t; x_0)$. Evaluating this integral, we obtain
\begin{equation}\label{Eq:WkSwim_S0}
    S_0(t; x_0) = 2\sum_{n=0}^\infty \frac{(-1)^n}{\lambda_{1,n}} \cos \left( \lambda_{1,n} \, x_0 \right)   e^{-\lambda_{1,n}^2 \, t }.
\end{equation}
To solve the $\mathcal{O}(Pe)$ problem \eqref{Eq:WkSwim_PesOrder_Prob}, 
we guess a solution of the form
\begin{align}\label{Eq:WeakSwim_Sol_n1_guess}
  n_1(x,t) = \sum_{n=0}^{\infty} \Big{[} \cos \left( \lambda_{1,n} \, x \right)p_n(t) +  \sin \left( \lambda_{2,n} x \right) q_n(t) \Big{]},
\end{align}
where $p_n(t)$ and $q_n(t)$ are functions to be determined. Note that \eqref{Eq:WeakSwim_Sol_n1_guess} satisfies the boundary condition \eqref{Eq:WkSwim_Pes_Order_BC} by construction. 
At $t=0$, we have $n_1(x,0) = 0$, which implies that $p_n(0)=q_n(0) = 0$ for all $n$. 
Integrating $n_1$ as given in \eqref{Eq:WeakSwim_Sol_n1_guess} over the domain, we obtain the survival probability at $\mathcal{O}(Pe_s)$, given by
\begin{align}
\label{Eq:S1-pn}
    S_1(t) = 2 \sum_{n=0}^\infty  \frac{(-1)^n}{\lambda_{1,n}} \; p_n(t).
\end{align}
Observe from the PDE for $n_1(x,t)$  in \eqref{Eq:WkSwim_PesOrder_Prob} that we need $f_0(x,t)$ to solve for $n_1$, since the equation for $n_1$ contains a derivative of $f_0$ with respect to $x$. Using a similar approach used to solve for $n_0$, we solve the $f_0(x,t)$ problem in \eqref{Eq:WkSwim_LD_Prob} to get
\begin{equation}\label{Eq:WeakSwim_Sol_f0}
\begin{split}
       f_0(x,t) = & (2\eta-1) \sum_{n=0}^{\infty} \Big{[} \cos \left( \lambda_{1,n} \, x_0 \right) \cos \left( \lambda_{1,n} \, x \right)e^{- \left( 2\beta + \lambda_{1,n}^2 \right)\, t }  + \sin(\lambda_{2,n} x_0) \sin \left( \lambda_{2,n} x \right) e^{-(2\beta + {\lambda_{2,n}^2)} \, t } \Big{]}.
  \end{split}
\end{equation}
Upon substituting \eqref{Eq:WeakSwim_Sol_n1_guess} and \eqref{Eq:WeakSwim_Sol_f0} into \eqref{Eq:WeakSwim_Pes_Order_n} and using the orthogonality properties of the sine and cosine functions, we derive ordinary differential equations (ODEs) for $p_n(t)$ and $q_n(t)$, given by
\begin{align}
  p_n'(t) +  \lambda_{1,n}^2 \, p_n(t) &= -(2\eta-1)\sum_{m=0}^\infty\ 2m\,\sin(\lambda_{2,m} \,x_0) A_{m,n} \,e^{-\big{(}2\beta+\lambda_{2,m}^2 \big{)} t}; \label{Eq:WeakSwim_pn_prob}\\[2ex]
    q_n'(t) +  \lambda_{2,n}^2 \,q_n(t) = & -(2\eta-1) \sum_{m=0}^\infty (2m+1) \cos\left( \lambda_{1,m} \,x_0\right)B_{m,n}\,e^{- \left( 2\beta +\lambda_{1,m}^2 \right)\, t }, \label{Eq:WeakSwim_qn_prob}
\end{align}
where $\lambda_{1,m} = (2m + 1)\pi/2$ and $\lambda_{2,m} = m\pi$ for $m \in \mathbb{Z}$, and $ A_{m,n}$ and $B_{m,n}$ for $n,m = 0, 1,2,\dots$ are defined as
\begin{align}\label{Eq:A_B_Integral}
        A_{m,n} &= \frac{\pi}{2} \int_{-1}^1\cos(\lambda_{2,m}x)  \cos(\lambda_{1,n} x) \;\text{d}x  \qquad \text{and} \qquad 
        B_{m,n} = \frac{\pi}{2} \int_{-1}^1\sin(\lambda_{1,m}x)  \sin(\lambda_{2,n} x) \;\text{d} x. 
\end{align}
Evaluating the integrals in \eqref{Eq:A_B_Integral} gives
\begin{align}\label{Eq:Anm_Bnm}
    A_{m,n} =  \frac{(-1)^{m+n}}{1+2(m+n)} + \frac{(-1)^{m-n}}{1-2(m-n)} \quad \text{and} \quad
    B_{m,n} = -\frac{(-1)^{m+n}}{1+2(m+n)} + \frac{(-1)^{m-n}}{1+2(m-n)}.
\end{align}
We impose the initial conditions $p_n(0)=0$ and $q_n(0) = 0$ for all $n = 0,1,2,\cdots$, on the ODEs in \eqref{Eq:WeakSwim_pn_prob} and \eqref{Eq:WeakSwim_qn_prob}. Using the method of integrating factor, we solve these ODE problems to obtain
\begin{align}
  p_n(t)  &= - \sum_{m=0}^\infty \frac{2m(2\eta-1) \sin(\lambda_{2,m}\, x_0) A_{m,n} }{\lambda_{1,n}^2 - (2\beta+\lambda_{2,m}^2)} \left[  e^{-\big{(} 2\beta+\lambda_{2,m}^2 \big{)}t}-e^{-  \lambda_{1,n}^2  t } \right]; \label{Eq:WeakSwim_pn_Sol2}\\[2ex]
    q_n(t)  &= -\sum_{m=0}^\infty  \frac{(2m+1)(2\eta-1)\cos(\lambda_{1,m} \,x_0)B_{m,n}}{\lambda_{2,n}^2-(2\beta + \lambda_{1,m}^2)}  \left[  e^{-(2\beta+\lambda_{1,m}^2)t}-e^{- \lambda_{2,n}^2 t } \right]. \label{Eq:WeakSwim_qn_Sol2}
\end{align}
Upon substituting \eqref{Eq:WeakSwim_pn_Sol2} into \eqref{Eq:S1-pn}, we obtain the survival probability at $\mathcal{O}(Pe)$ as 
\begin{subequations}\label{Eq:WkSwim_S1}
    \begin{align}\label{Eq:WkSwim_S1_A}
    S_1(t; x_0) = \sum_{n=0}^\infty    \sum_{m=0}^\infty \Psi_{n}\,\Phi_{n,m} \left(  e^{-\big{(}2\beta+\lambda_{2,m}^2 \big{)}t}-e^{-  \lambda_{1,n}^2  t } \right),
\end{align}
where 
\begin{align}\label{Eq:Psi_Phi}
   \Psi_{n} =  \frac{2(-1)^n}{\lambda_{1,n}}\quad \text{and} \quad \Phi_{n,m} =  \frac{2m(1-2\eta) \sin(\lambda_{2,m}\, x_0) A_{m,n} }{\lambda_{1,n}^2 - (2\beta+\lambda_{2,m}^2)} .
\end{align}
\end{subequations}
From \eqref{Eq:WkSwim_S0} and \eqref{Eq:WkSwim_S1}, we construct a two-term asymptotic expansion for the survival probability of the ABPs in $\Omega_u$ in the weak-swimming ($Pe \ll 1$) regime, given by
\begin{equation}\label{Eq:WkSwim_S-two_term}
\begin{split}
    S(t;x_0) &=  \sum_{n=0}^\infty \Psi_{n}  \cos \left(\lambda_{1,n} \, x_0 \right)   e^{-\lambda_{1,n}^2 \, t }  + Pe \sum_{n=0}^\infty  \Psi_{n}  \sum_{m=0}^\infty \,\Phi_{n,m} \left(  e^{-(2\beta+\lambda_{2,m}^2)\,t}-e^{-  \lambda_{1,n}^2 \, t } \right) +\mathcal{O}(Pe^2),
\end{split}    
\end{equation}
We remark that the leading-order term in the survival probability in \eqref{Eq:WkSwim_S-two_term} corresponds to the survival probability of a passive Brownian particle in $\Omega_u$, and the first effect of swimming on the survival probability comes at the  $\mathcal{O}(Pe)$ term. Integrating \eqref{Eq:WkSwim_S-two_term}, we construct a two-term asymptotic expansion for the MFPT (cf. \eqref{Eq:MFPT_expand}), $\mu(x_0) = \mu_0 + Pe \mu_1 +\mathcal{O}(Pe^2)$, 
where 
\begin{equation}
\label{Eq:WkSwm:2Term_MFPT}
    \mu_0 = 2\sum_{n=0}^\infty \frac{(-1)^n}{\lambda_{1,n}^3} \cos \left( \lambda_{1,n}  \, x_0 \right) \quad \mathrm{and}\quad \mu_1 = \sum_{n=0}^\infty  \frac{2(-1)^n}{\lambda_{1,n}^3} \sum_{m=0}^\infty  \frac{2m(1-2\eta)A_{m,n}}{2\beta + \lambda_{2,m}^2}\sin(\lambda_{2,m}\, x_0).
\end{equation}
In the preceding equation, $\lambda_{1,n} =  (2n + 1)\pi/2$, $\lambda_{2,n} = n\pi$ and $A_{m,n}$ is as defined in \eqref{Eq:Anm_Bnm}.

We observe that when the particles start from the midpoint of the domain  $\Omega_u$ ($x_0 = 0$), the $\Order{Pe}$ term in the MFPT expansion vanishes due to symmetry, i.e., $\mu_1 =0$. In this case, the first effect of swimming on the MFPT comes at $\Order{Pe^2}$, as we shall show later. We also note that regardless of the starting position $x_0$,  $\mu_1=0$ when $\eta=1/2$. That is, the $\Order{Pe}$ MFPT vanishes if we start with the particles having an equal probability of pointing to the positive and negative sides (unbiased starting orientations). To see this, first notice that $f_0 =0$ if $\eta=1/2$; as a result, from \eqref{Eq:WkSwim_PesOrder_Prob} we then obtain  $n_1= 0$ and accordingly $\mu_1 =0$.

As noted above, when $x_0 = 0$ or $\eta=1/2$, the $\Order{Pe}$ MFPT vanishes and we need to continue the asymptotic analysis to higher order. To this end, we consider the $\mathcal{O}(Pe^2)$ problem given in \eqref{Eq:WkSwim_Pes2Order_Prob}. To solve this problem, we need to first determine $f_1(x,t)$. We consider the PDE in \eqref{Eq:WeakSwim_Pes_Order_f} and guess a solution of the form
\begin{equation}\label{Eq:f1_Sol_guess}
    f_1(x,t) = \sum_{n=0}^{\infty} \cos(\lambda_{1,n}\,x)\; u_{n}(t) + \sin(\lambda_{2,n}\,x)\; w_{n}(t),
\end{equation}
where we recall that $\lambda_{1,n} = (2n + 1)\pi/2$ and $\lambda_{2,n} = n\pi$ for $n \in \mathbb{Z}$, and $u_{n}(t)$ and $w_{n}(t)$ for $n = 0, 1,2,\dots$ are functions to be determined. Upon substituting $n_0(x,t)$ as given in \eqref{Eq:WkSwim_n0_Sol} and $f_1(x,t)$ given in \eqref{Eq:f1_Sol_guess} into \eqref{Eq:WeakSwim_Pes_Order_f}, we use the orthogonality properties of sine and cosine on $[-1, 1]$ to derive  ODEs for $u_{n}(t)$ and $w_{n}(t)$, given by
\begin{equation}\label{Eq:WkSwim_Un_Wn_ODE}
\begin{split}
   \dot{u}_{n} + \chi_{1,n}^2\,u_{n} &= - \frac{2}{\pi} \sum_{m=0}^{\infty} \lambda_{2,m}\, A_{m,n}\,\sin(\lambda_{2,m}\,x_0) \, e^{-\lambda_{2,m}^2 \, t}, \qquad  u_{n}(0) = 0\; ; \\[1ex]
  \dot{w}_{n} + \chi_{2,n}^2 \,w_{n} &=  \frac{2}{\pi} \sum_{m=0}^{\infty} \lambda_{1,m}\, B_{m,n}\,\cos(\lambda_{1,m}\,x_0) \, e^{-\lambda_{1,m}^2 \, t}, \qquad  w_{n}(0) = 0\;,
\end{split}    
\end{equation}
where  $\chi_{1,n}^2 = \lambda_{1,n}^2 + 2\beta$, $\chi_{2,n}^2 = \lambda_{2,n}^2 + 2\beta$, and $ A_{m,n}$ and $B_{m,n}$ for $n,m = 0, 1,2,\dots$ are as defined in \eqref{Eq:Anm_Bnm}. Solving the ODEs in \eqref{Eq:WkSwim_Un_Wn_ODE}, we obtain
\begin{equation}\label{Eq:WkSwim_Un_Wn_Sol}
\begin{split}
   u_{n} &= - \frac{2}{\pi} \sum_{m=0}^{\infty} \lambda_{2,m}\,  \frac{A_{m,n}\,\sin(\lambda_{2,m}\,x_0)}{(\chi_{1,n}^2 - \lambda_{2,m}^2 )} \,\left( e^{-\lambda_{2,m}^2 \, t}  - e^{-\chi_{1,n}^2 \, t} \right),\\[1ex]
  w_{n} &=  \frac{2}{\pi} \sum_{m=0}^{\infty} \lambda_{1,m}\, \frac{B_{m,n}\,\cos(\lambda_{1,m}\,x_0)}{(\chi_{2,n}^2 - \lambda_{1,m}^2 )} \, \left( e^{-\lambda_{1,m}^2 \, t} - e^{-\chi_{2,n}^2 \, t} \right),
\end{split}    
\end{equation}
We substitute
\eqref{Eq:WkSwim_Un_Wn_Sol} into \eqref{Eq:f1_Sol_guess} to construct the solution for $f_1(x,t)$, given by
\begin{equation}\label{Eq:f1_Sol}
\begin{split}
    f_1(x,t) &= - \frac{2}{\pi}\sum_{n=0}^{\infty} \sum_{m=0}^{\infty} \left[ \cos(\lambda_{1,n}\,x)\;   \lambda_{2,m}\,  \frac{A_{m,n}\,\sin(\lambda_{2,m}\,x_0)}{(\chi_{1,n}^2 - \lambda_{2,m}^2 )} \,\left( e^{-\lambda_{2,m}^2 \, t}  - e^{-\chi_{1,n}^2 \, t} \right) \right. \\
    & \qquad \qquad \left. - \; \sin(\lambda_{2,n}\,x)\, \lambda_{1,m}\, \frac{B_{m,n}\,\cos(\lambda_{1,m}\,x_0)}{(\chi_{2,n}^2 - \lambda_{1,m}^2 )} \, \left( e^{-\lambda_{1,m}^2 \, t} - e^{-\chi_{2,n}^2 \, t} \right) \right].
\end{split}      
\end{equation}

Now, we consider the PDE for $n_2(x,t)$ given in \eqref{Eq:WkSwim_Pes2Order_Prob} and guess a solution of the form
\begin{equation}\label{Eq:n2_Sol_guess}
    n_2(x,t) = \sum_{n=0}^{\infty} \cos(\lambda_{1,n}\,x)\; g_{n}(t) + \sin(\lambda_{2,n}\,x)\; h_{n}(t),
\end{equation}
where  $\lambda_{1,n} = (2n + 1)\pi/2$ and $\lambda_{2,n} = n\pi$ for $n \in \mathbb{Z}$, and $g_{n}(t)$ and $h_{n}(t)$ for $n = 0, 1,2,\dots$ are functions to be determined. Substituting \eqref{Eq:f1_Sol} and \eqref{Eq:n2_Sol_guess} into \eqref{Eq:WeakSwim_Pes2_Order_n} and using the orthogonality properties of sine and cosine on $[-1, 1]$, we construct ODE problems for $g_{n}(t)$ and $h_{n}(t)$, given by
\begin{equation}\label{Eq:WkSwim_gn_hn_ODE}
\begin{split}
   \dot{g}_{n} + \lambda_{1,n}^2\,g_{n} &= - \frac{4}{\pi^2} \sum_{m=0}^{\infty} \lambda_{2,m}\, A_{m,n}\,\sum_{k=0}^{\infty} \lambda_{1,k}\, \frac{B_{k,m}\,\cos(\lambda_{1,k}\,x_0)}{(\chi_{2,m}^2 - \lambda_{1,k}^2 )} \, \left( e^{-\lambda_{1,k}^2 \, t} - e^{-\chi_{2,m}^2 \, t} \right), \qquad  g_{n}(0) = 0\; ; \\[1ex]
  \dot{h}_{n} + \lambda_{2,n}^2 \,h_{n} &=  -\frac{4}{\pi^2} \sum_{m=0}^{\infty} \lambda_{1,m}\, B_{m,n}\, \sum_{k=0}^{\infty} \lambda_{2,k}\,  \frac{A_{k,m}\,\sin(\lambda_{2,k}\,x_0)}{(\chi_{1,m}^2 - \lambda_{2,k}^2 )} \,\left( e^{-\lambda_{2,k}^2 \, t}  - e^{-\chi_{1,m}^2 \, t} \right), \qquad  h_{n}(0) = 0\;,
\end{split}    
\end{equation}
where $ A_{m,n}$ and $B_{m,n}$ for $n,m = 0, 1,2,\dots$ are as defined in \eqref{Eq:Anm_Bnm}. We solve the ODEs in \eqref{Eq:WkSwim_gn_hn_ODE} to obtain
\begin{equation}\label{Eq:WkSwim_gn_hn_Sol}
\begin{split}
   g_{n}(t) &= - \frac{4}{\pi^2} \sum_{m=0}^{\infty} \lambda_{2,m}\, A_{m,n}\,\sum_{k=0}^{\infty} \frac{\lambda_{1,k}\, B_{k,m}\,\cos(\lambda_{1,k}\,x_0)}{(\chi_{2,m}^2 - \lambda_{1,k}^2 )} \left[
  \frac{ \left( e^{-\lambda_{1,k}^2 \, t} - e^{- \lambda_{1,n}^2 \, t} \right)}{(\lambda_{1,n}^2 - \lambda_{1,k}^2 )} - 
  \frac{ \left( e^{-\chi_{2,m}^2 \, t} - e^{- \lambda_{1,n}^2 \, t} \right)}{(\lambda_{1,n}^2 - \chi_{2,m}^2 )}  \right],  \\[1ex]
h_{n}(t) &=  -\frac{4}{\pi^2} \sum_{m=0}^{\infty} \lambda_{1,m}\, B_{m,n}\, \sum_{k=0}^{\infty}   \frac{\lambda_{2,k}\,A_{k,m}\,\sin(\lambda_{2,k}\,x_0)}{(\chi_{1,m}^2 - \lambda_{2,k}^2 )} \,\left[
  \frac{ \left( e^{-\lambda_{2,k}^2 \, t} - e^{- \lambda_{2,n}^2 \, t} \right)}{(\lambda_{2,n}^2 - \lambda_{2,k}^2 )} - 
  \frac{ \left( e^{-\chi_{1,m}^2 \, t} - e^{- \lambda_{2,n}^2 \, t} \right)}{(\lambda_{2,n}^2 - \chi_{1,m}^2 )}  \right].
\end{split}    
\end{equation}
Therefore, the $\mathcal{O}(Pe^2)$ density of particles $n_2(x,t) $ can be derived explicitly be substituting  the functions $g_n(t)$ and $h_n(t)$ as given in \eqref{Eq:WkSwim_gn_hn_Sol} into \eqref{Eq:n2_Sol_guess}. Integrating the resulting expression with respect to $x$ over the domain $[-1,1]$, we obtain the survival probability
\begin{equation}\label{Eq:WkSwim_Surv2_Sol}
\begin{split}
   S_2(t) &= - \frac{8}{\pi^2} \sum_{n=0}^{\infty} \frac{(-1)^n}{\lambda_{1,n}} \sum_{m=0}^{\infty} \lambda_{2,m}\, A_{m,n}\,\sum_{k=0}^{\infty} \frac{\lambda_{1,k}\, B_{k,m}\,\cos(\lambda_{1,k}\,x_0)}{(\chi_{2,m}^2 - \lambda_{1,k}^2 )} \left[
  \frac{ \left( e^{-\lambda_{1,k}^2 \, t} - e^{- \lambda_{1,n}^2 \, t} \right)}{(\lambda_{1,n}^2 - \lambda_{1,k}^2 )} - 
  \frac{ \left( e^{-\chi_{2,m}^2 \, t} - e^{- \lambda_{1,n}^2 \, t} \right)}{(\lambda_{1,n}^2 - \chi_{2,m}^2 )}  \right],  
\end{split}    
\end{equation}
Note that $h_n(t)$ vanishes from the survival probability 
since $\int_{-1}^{1}  \sin(\lambda_{2,n}\,x) \; \text{d}x = 0$. Lastly, we integrate $S_2(t)$ with respect to $t$ to obtain the MFPT at $\mathcal{O}(Pe^2)$, given by 
\begin{equation}\label{Eq:WkSwim_MFPT2}
\begin{split}
   \mu_2(x_0) &=  \frac{8}{\pi^2} \sum_{n=0}^{\infty} \frac{(-1)^{n+1}}{\lambda_{1,n}^3} \sum_{m=0}^{\infty} \frac{\lambda_{2,m}}{\chi_{2,m}^2}\, A_{m,n}\,\sum_{k=0}^{\infty} \frac{ B_{k,m}\,\cos(\lambda_{1,k}\,x_0)}{  \lambda_{1,k}} .
\end{split}     
\end{equation}
We combine  \eqref{Eq:WkSwm:2Term_MFPT} and \eqref{Eq:WkSwim_MFPT2}  to construct a three-term asymptotic expansion for the MFPT,  given by 
 \begin{equation}\label{Eq:WkSwm:3Term_MFPT}
\begin{split}
 \mu(x_0) &= 2\sum_{n=0}^\infty \frac{(-1)^n}{\lambda_{1,n}^3} \cos \left( \lambda_{1,n}  \, x_0 \right)  + Pe \sum_{n=0}^\infty  \frac{2(-1)^n}{\lambda_{1,n}^3} \sum_{m=0}^\infty  \frac{2m(1-2\eta)A_{m,n}}{2\beta + \lambda_{2,m}^2}\sin(\lambda_{2,m}\, x_0) \\[1ex]
 & \qquad + Pe^2 \; \sum_{n=0}^{\infty} \frac{8}{\pi^2}\frac{(-1)^{n+1}}{\lambda_{1,n}^3} \sum_{m=0}^{\infty} \frac{\lambda_{2,m}}{\chi_{2,m}^2}\, A_{m,n}\,\sum_{k=0}^{\infty} \frac{ B_{k,m}\,\cos(\lambda_{1,k}\,x_0)}{\lambda_{1,k} }  + \mathcal{O}(Pe^3).
\end{split}
\end{equation}

\begin{figure}[!h]
  \centering
  \includegraphics[scale=1.07]{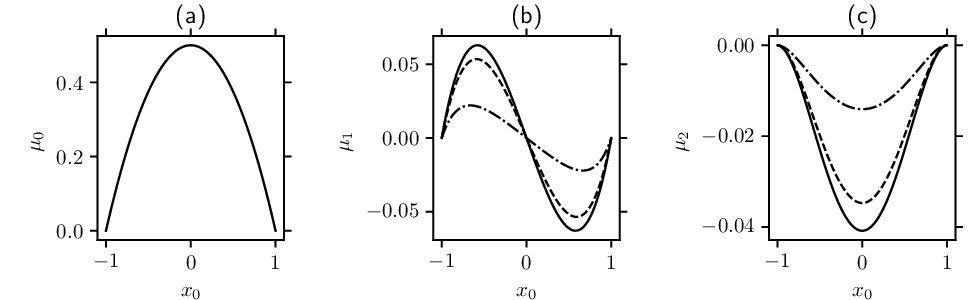}
  \caption{\textbf{Contribution of the asymptotic terms to the MFPT.} Plots of each term in the asymptotic expansion \eqref{Eq:WkSwm:3Term_MFPT} as a function of the particles' starting position ($x_0$) for different values of $\beta$: $\beta = 0.1$ (solid line), $\beta = 1$ (dashed line), and $\beta = 10$ (dash-dot line). (a) Leading-order MFPT $(\mu_0)$, (b) $\mathcal{O}(Pe)$ MFPT ($\mu_1)$, and (c) $\mathcal{O}(Pe^2)$ MFPT ($\mu_2$). Results were computed with 100 terms of each series in \eqref{Eq:WkSwm:3Term_MFPT} and $\eta=1$.}
\label{Fig:MFPT_WkSwim_MuVsX0}
\end{figure}

Figure~\ref{Fig:MFPT_WkSwim_MuVsX0} shows the plot for each term in the asymptotic expansion \eqref{Eq:WkSwm:3Term_MFPT} for different values of $\beta$: $\beta = 0.1$ (solid line), $\beta = 1$ (dashed line), and $\beta = 10$ (dash-dot line),  as a function of the  starting position, $x_0$. The leading-order MFPT is shown on the left, the   $\mathcal{O}(Pe)$ term in the middle, and the  $\mathcal{O}(Pe^2)$ term on the right.
For these results, truncated sums are computed using 100 terms of each series in \eqref{Eq:WkSwm:3Term_MFPT} and $\eta = 1$, which corresponds to all the particles starting with positive orientation, i.e., pointing to the right. The leading-order MFPT is the (passive) Brownian result, and thus does not depend on $\beta$ (cf. \eqref{Eq:WkSwim_LD_Prob}). The Brownian MFPT obtains its maximum when particles are started at the center-point of the 1-D domain and decreases as the starting point moves closer to the boundaries. The Brownian MFPT profile as a function of $x_0$ is symmetric. The first effect of swimming on the MFPT is obtained at  $\mathcal{O}(Pe)$. We observe that the contribution of the $\mathcal{O}(Pe)$ correction term to the MFPT depends on the starting position with the least contribution when the particles start close to the two ends of the domain or in the middle of the domain (at the origin). More importantly, at  $\mathcal{O}(Pe)$, the swimming motion can either increase or decrease the MFPT depending on the starting position. Because initially all particles started pointing to the right ($\eta=1$), the MFPT is decreased when $x_0 \in (0,1)$ and increases when $x_0 \in (-1,0)$. The reduction in the MFPT for $x_0 \in (0,1)$ is because the particles can easily swim to the right boundary and exit the domain at $x=1$ since they started with positive orientation. As can be seen in Figure~\ref{Fig:MFPT_WkSwim_MuVsX0}b, there is an optimal value of the starting position that maximizes the reduction of the MFPT at $\mathcal{O}(Pe)$. The $\mathcal{O}(Pe^2)$ correction term decreases the MFPT regardless of the starting position. 

We also observe from the results in Figure~\ref{Fig:MFPT_WkSwim_MuVsX0}  that the effects of the  $\mathcal{O}(Pe)$ and $\mathcal{O}(Pe^2)$ correction terms on the MFPT decrease as $\beta$ increases.
Recall that $\beta= R^2 /(\tau D_T)$, where $R$ is the length scale of the domain, $D_T$ is the dimensional diffusion rate of the particles, and $1/\tau$ is the tumbling rate of the particles (the rate at which the particles switch their orientation).
For fixed values of $R$ and $D_T$, an increase in $\beta$ corresponds to an increase in the tumbling rate of the particles. When the particles  tumble fast, they lose their persistent swim motion and behave more like that of Brownian particles. As a result, the effect of the correction terms in \eqref{Eq:WkSwm:3Term_MFPT} decreases as $\beta$ increases.

\begin{figure}[!h]
  \centering
  \includegraphics[scale=0.9]{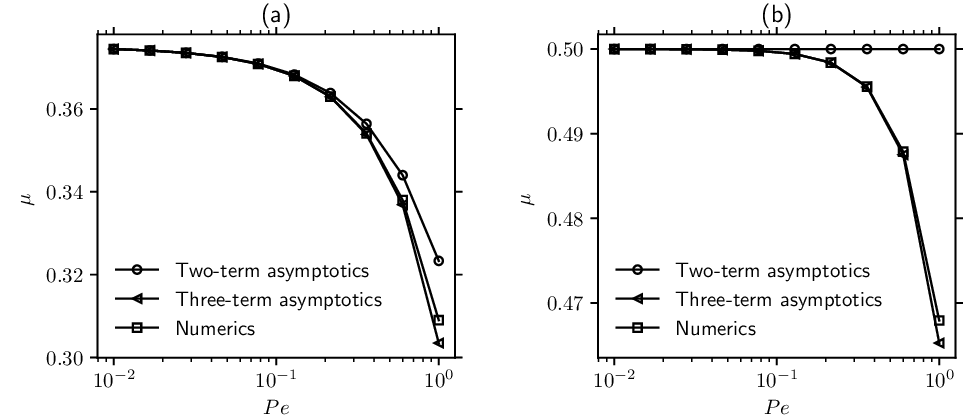}
  \caption{\textbf{Asymptotic and numerical MFPT (low tumbling rate).} Comparisons of the asymptotic approximation \eqref{Eq:WkSwm:3Term_MFPT} and the numerical solution of the full PDE \eqref{Eq:Dimless_Model_n}.
    (a) Particles' starting position, $x_0 = 0.5$ and $\eta$ = 1 (all the particles were initially pointing to the right). (b) Particles' starting position, $x_0 = 0$ and $\eta = 0.5$ (half of the particles were initially pointing in opposite directions). Results were computed with 100 terms of each series in \eqref{Eq:WkSwm:3Term_MFPT} and $\beta=1$.}
  \label{Fig:MFPT_WkSwim_Asyp_vs_Num}
\end{figure}

The results in Figure~\ref{Fig:MFPT_WkSwim_Asyp_vs_Num} show comparisons of the analytical MFPT ($\mu$) in \eqref{Eq:WkSwm:3Term_MFPT} and the numerical solution of the full PDE \eqref{Eq:Dimless_Model_n} as a function of the P\'eclect number ($Pe$) for (a) $x_0 = 0.5$ and $\eta = 1$ and (b) $x_0 = 0$ and $\eta = 0.5$. For each plot in this figure, we show the two- and three-term approximation of the MFPT \eqref{Eq:WkSwm:3Term_MFPT} and the numerical MFPT. 
These results show that the MFPT decreases as $Pe$ increases. Recall that $Pe = (v_s R)/D_T$, which implies that $Pe$ is directly proportional to the dimensional swimming speed of the particles ($v_s$). 
For fixed values of $R$ and $D_T$, increasing $Pe$ corresponds to increasing the swimming speed of the particles. 
Therefore, as $Pe$ increases, the particles swim faster, making it easier to exit the domain and decrease the MFPT, as shown in the results.
The asymptotic approximation agrees well with the numerical solution for small values of $Pe$. 
However, the two solutions deviate as $Pe$ increases. As expected, the three-term asymptotic MFPT agrees better with the numerical solution for higher  values of $Pe$ than the two-term expansion.
In Figure~\ref{Fig:MFPT_WkSwim_Asyp_vs_Num}a, where the particles start from $x_0 = 0.5$ and point to the right, the two-term asymptotic expansion agrees with the numerical solution up to  $Pe \approx  0.2$, after which it overestimates the MFPT as $Pe$ continues to increase. On the other hand, the three-term asymptotic approximation agrees with the numerical solution up to $Pe \approx 0.7$ and then underestimates the MFPT as $Pe$ increases.

In Figure~\ref{Fig:MFPT_WkSwim_Asyp_vs_Num}b, the particles started at the origin ($x_0 = 0$), with half of them pointing to the right while the other half is pointing to the left, i.e., $\eta = 1/2$.
As mentioned earlier, the $\mathcal{O}(Pe)$ correction term in the MFPT expansion \eqref{Eq:WkSwm:3Term_MFPT} varnishes when the particles start from the origin.
As a result of this, the two-term expansion of the MFPT reduces to the leading-order MFPT, which is independent of the particle's swimming speed or the P\'eclet number. 
This corresponds to the horizontal line in Figure~\ref{Fig:MFPT_WkSwim_Asyp_vs_Num}b. 
Similar to the results in Figure~\ref{Fig:MFPT_WkSwim_Asyp_vs_Num}a, the three-term MFPT approximation underestimates the MFPT as $Pe$ increases. We used 100 terms of the series in \eqref{Eq:WkSwm:3Term_MFPT} and $\beta = 1$ for the results in this figure.
\begin{figure}[!h]
  \centering
  \includegraphics[scale=0.9]{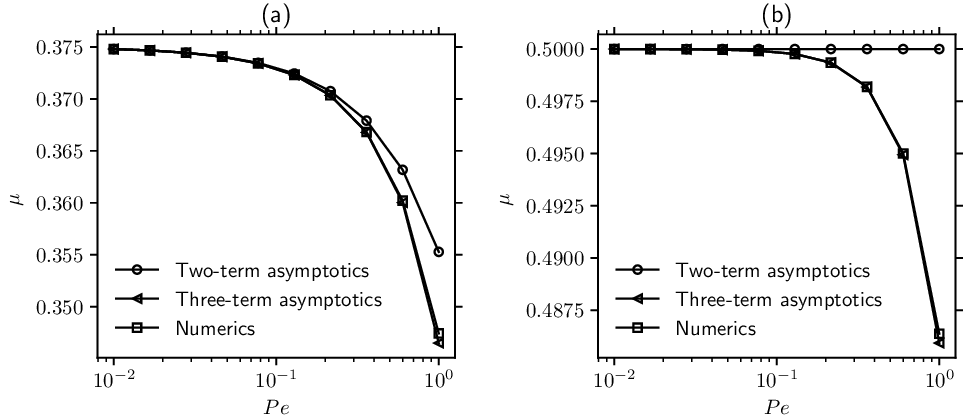}
  \caption{\textbf{Asymptotic and numerical MFPT (high tumbling rate).} Comparisons of the asymptotic approximation \eqref{Eq:WkSwm:3Term_MFPT} and the numerical solution of the full PDE \eqref{Eq:Dimless_Model_n}.
   (a) Particles' starting position, $x_0 = 0.5$ and $\eta$ = 1 (all the particles were initially pointing to the right). (b) Particles' starting position, $x_0 = 0$ and $\eta = 0.5$ (half of the particles were initially pointing in opposite directions). Results were computed with 100 terms of each series in \eqref{Eq:WkSwm:3Term_MFPT} and $\beta=10$. }
\label{Fig:MFPT_WkSwim_Asyp_vs_Num_beta}  
\end{figure}

In Figure~\ref{Fig:MFPT_WkSwim_Asyp_vs_Num_beta}, we repeat the comparisons shown in Figure~\ref{Fig:MFPT_WkSwim_Asyp_vs_Num} but for $\beta = 10$, i.e., a higher tumbling rate. Here,  the asymptotic MFPT results agree more with the numerical result than those in Figure~\ref{Fig:MFPT_WkSwim_Asyp_vs_Num}. In addition, the three-term expansion agrees more with the numerical result than the two-term expansion. 

\begin{figure}[!h]
  \centering
  \includegraphics[scale=1.]{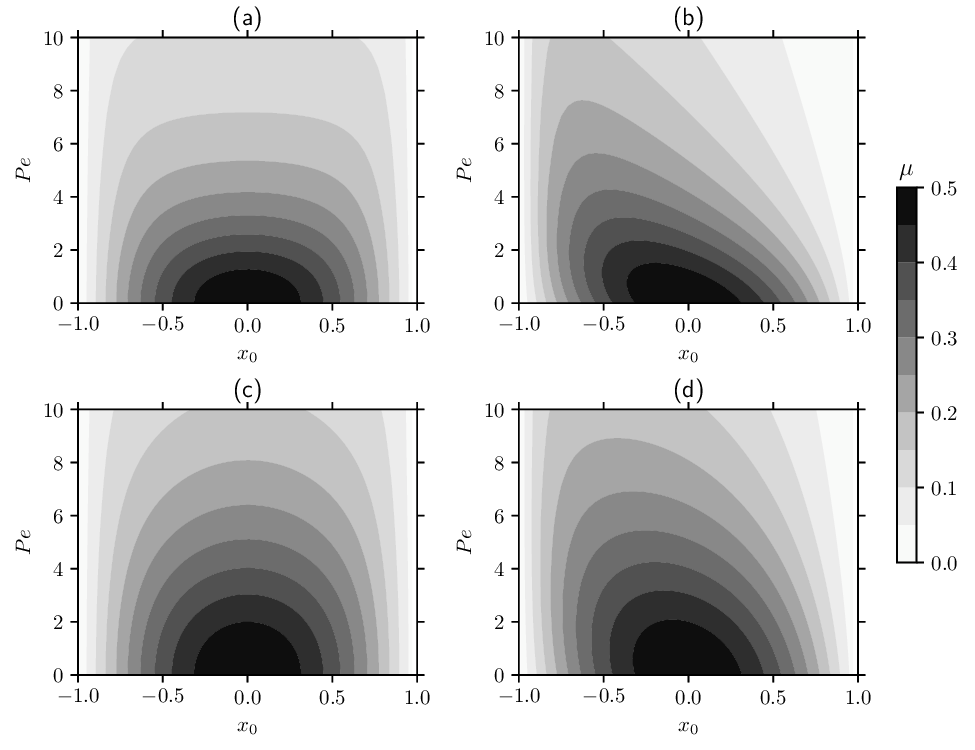}
  \caption{\textbf{MFPT computed from the full PDE.} Contour plots of MFPT ($\mu$) computed from numerical simulations of the full PDE \eqref{Eq:Dimless_Model_n} for different starting position ($x_0$) and P\'eclet number ($Pe$). (a) $\beta=1$ and $\eta=0.5$, (b) $\beta=1$ and $\eta=1$,  (c) $\beta=10$ and $\eta=0.5$,  and (d) $\beta=10$ and  $\eta=1$. Note: $\eta = 1$ implies that all the particles were initially pointing to the right, and $\eta= 0.5$ means half of the particles were initially pointing to the right while the other half were pointing to the left.}
  \label{Fig:MFPT_numerics_contour}
\end{figure}

Next, we present contour plots of the MFPT \eqref{Eq:WkSwm:3Term_MFPT} as a function of the particles' starting position ($x_0$) and P\'eclet number ($Pe$) for different scenarios:  (a) $\beta=1$ and $\eta=0.5$, (b) $\beta=1$ and $\eta=1$,  (c) $\beta=10$ and $\eta=0.5$,  and (d) $\beta=10$ and  $\eta=1$.  For Figures~\ref{Fig:MFPT_numerics_contour}a and \ref{Fig:MFPT_numerics_contour}c, where half of the particles started by facing opposite directions ($\eta = 0.5$), the MFPT is symmetric about the center of the domain. In this symmetric case, the maximum in MFPT is achieved when the particles start from the center ($x_0=0$) of the domain. The variations in the MFPT is weak  for $x_0$ values close to the origin when $\beta = 1$, and this range of $x_0$ increases as $Pe$ increases (Figure~\ref{Fig:MFPT_numerics_contour}a). However, the MFPT decreases more rapidly as $x_0$ shifts away from the origin, towards the two domain boundaries (Figure~\ref{Fig:MFPT_numerics_contour}c) when $\beta$ is increased to $\beta = 10$.

In Figures~\ref{Fig:MFPT_numerics_contour}b and \ref{Fig:MFPT_numerics_contour}d, where all the particles started by pointing to the right, the MFPT is not symmetric about the origin. Rather it is  lower for $x_0 \in (0, 1)$, compared to the second half of the domain, $(-1, 0)$. As the particles' swimming speed increases, the effect of the initial orientation of the particles on the MFPT increases, making the lack of symmetry in the MFPT become more apparent. Since the particles started by pointing to the right, they can easily exit the domain when they start from the right half of the domain, compared to the left half of the domain. 
 Although the results in Figures~\ref{Fig:MFPT_numerics_contour}b and \ref{Fig:MFPT_numerics_contour}d were computed for particles that started by pointing to the right, analogous results are expected when all the particles start by pointing to the left. In that case, one would expect the MFPT to be lower for the particles starting on the left half of the domain than those starting from the right half.

\section{Discussion}\label{sec:Discuss}

The first passage time distribution and its mean have been studied extensively for passive Brownian particles \cite{pillay2010asymptotic, cheviakov2010asymptotic, iyaniwura2021optimization, 
lindsay2017first, iyaniwura2021simulation} and applied to problems in physical and biological sciences. However, for many biological systems and processes, the elements are often active and exhibit self-propulsion in addition to Brownian motion. Despite the commonality of active Brownian motion,  not much has been done in terms of developing a first passage time framework for these processes until recently \cite{biswas2023escape, bressloff2023close, bressloff2023trapping, di2023active}.
Here, we used asymptotic analysis to study the MFPT for ABPs in a 1-D domain with absorbing boundaries at the two ends in a weak-swimming limit.
We construct a three-term asymptotic approximation for the  MFPT  in terms of the P\'eclet number, which is directly proportional to the swimming speed of the particles. Numerical simulations of the full PDE were used to validate our analytical approximation.

The leading-order term in our asymptotic approximation of the MFPT for ABPs represents the MFPT for Brownian particles, and the effect of swimming comes in at the $\mathcal{O}(Pe)$ correction term onward. However, our analysis shows that the $\mathcal{O}(Pe)$ term vanishes when the particles start at the origin ($x_0$) or when half of them start by  facing opposite directions ($\eta=1/2$). Note that the $\mathcal{O}(Pe)$ density of ABPs ($n_1$) is driven by $f_0$ (see \eqref{Eq:WkSwim_PesOrder_Prob}). When  $\eta=1/2$, $f_0$ vanishes, and as a result, we have $n_1=0$ and $\mu_1=0$. In these cases, the correction to the leading-order MFPT comes in at the $\mathcal{O}(Pe^2)$ term.
The contribution of the $\mathcal{O}(Pe)$ term to the MFPT is asymmetric about the center of the domain.
On the other hand, the $\mathcal{O}(Pe^2)$ is symmetric about the center of the domain, with the minimum contributions occurring when the particles start close to the two boundaries of the domain and the maximum contribution when they start at the origin. 
The contributions from both the $\mathcal{O}(Pe)$ and $\mathcal{O}(Pe^2)$ terms decrease as the tumbling rate of the particles increases (Figure~\ref{Fig:MFPT_WkSwim_MuVsX0}). As the tumbling rate of the particles increases, the motion of the ABPs becomes more diffusive, and the effect of swimming is minimized. As a result, the effect of the correction terms also decreases.

Although our analytical approximation of the MFPT is valid for small swimming speeds, it still agrees with the MFPT computed numerically by solving the full PDE model for reasonable values of $Pe$. As the swimming speed increases, the analytical MFPT deviates from the numerical MFPT.
In this study, we assume that the diffusion rate of the particles is constant. As expected,  the time it takes for the particles to exit the domain decreases as their swimming speed increases.
On the other hand, the MFPT of the active particles approaches that of Brownian particles as the swimming speed decreases.
As the tumbling rate of the particles increases, their motion becomes more diffusive, and as a result, the MFPT approaches that of Brownian particles. When the particles start with half of them pointing to the right and the other half pointing to the left, the MFPT is symmetric about the center of the domain. However, when all the particles start by pointing in a specific direction, the MFPT is skewed in that direction (Figure~\ref{Fig:MFPT_numerics_contour}). This helps to account for the effect of the particles' swimming and their initial orientation.
When all the particles start in the right half of the domain and point to the right, their MFPT is lower than those that started  by pointing in the same direction but from the left half of the domain. The higher MFPT for these particles is due to their initial orientation. As the particles' swimming speed decreases, this effect is minimized.

Our work provides some contributions to the growing literature on the MFPT of active Brownian particles. In particular, the weak-swimming asymptotic analysis allows us to reveal the first effect of swimming on the MFPT in contrast to that of passive particles. We have shown that both the starting position and orientation play important roles in the behavior of the MFPT of active particles. An interesting and straightforward extension to this work is to consider different boundary conditions. For example, one could consider a reflecting boundary condition at one end of the domain while the other end has an absorbing boundary condition. In this case, the particle can only escape the domain through the boundary that is absorbing. The mean first passage time, in this case, may not be trivial, especially when the particles start by pointing in the direction of the reflecting boundary. It would also be worthwhile to consider a similar analysis for a domain with partially absorbing boundaries or stochastic switching boundaries, such as those used in \cite{mercado2021first}. Future work on this topic should also include extending the problem to 2-D and 3-D domains.

\section*{Funding statement}

The authors received no financial support for the research, authorship, and/or publication of this article. 


\section*{Declaration of Competing Interest}

The authors declare that there is no competing interest.

\bibliographystyle{unsrtnat}
\bibliography{References}

\newpage

\renewcommand{\thesection}{S\arabic{section}}
\renewcommand{\thefigure}{S\arabic{figure}}
\renewcommand{\thetable}{S\arabic{table}}
  \setcounter{section}{0}
  \setcounter{figure}{0}
  \setcounter{table}{0}
\appendixpageoff
\appendixtitleoff
\renewcommand{\appendixtocname}{Supplementary material}

\clearpage
\pagenumbering{arabic}
\renewcommand*{\thepage}{S-\arabic{page}}

\end{document}